\documentclass[12pt,leqno]{amsart} 
\newtheorem{thm}{Theorem}[section]
\newtheorem{defi}{Definition}[section]

\newtheorem{cor}{Corollary}[section]
\newtheorem{pr}{Proposition}[section]

\newcommand{\be}{\begin{equation}}
\newcommand{\ee}{\end{equation}}
\newcommand{\bea}{\begin{eqnarray}}
\newcommand{\eea}{\end{eqnarray}}
\newcommand{\beb}{\begin{eqnarray*}}
\newcommand{\eeb}{\end{eqnarray*}}
\usepackage{amssymb,amsfonts,amsthm,setspace,indentfirst,mathrsfs}
\usepackage{minibox}
\numberwithin{equation}{section}
\begin{document}
%
\title[On super generalized recurrent manifolds]{\bf{On super generalized recurrent manifolds}}
\author[Absos Ali Shaikh, Indranil Roy and Haradhan Kundu]{Absos Ali Shaikh, Indranil Roy and Haradhan Kundu}

\address{\noindent\newline Department of Mathematics,\newline The University of 
Burdwan, Golapbag,\newline Burdwan-713104,\newline West Bengal, India}
\email{aask2003@yahoo.co.in, aashaikh@math.buruniv.ac.in}
\email{royindranil1@gmail.com}
\email{kundu.haradhan@gmail.com}
\dedicatory{}
\begin{abstract}
To generalize the notion of recurrent manifold, there are various recurrent like conditions in the literature. In this paper we present a recurrent like structure, namely, \textit{super generalized recurrent manifold}, which generalizes both the hyper generalized recurrent manifold and weakly generalized recurrent manifold. The main object of the present paper is to study the geometric properties of super generalized recurrent manifold. Finally to ensure the existence of such structure we present a proper example by a suitable metric.
\end{abstract}
%
\maketitle
%

\section{\bf{Introduction}}
 Let $M$ be a connected semi-Riemannian manifold equipped with a semi-Riemannian metric $g$. Let $\nabla$, $R$, $S$ and
 $\kappa$ be respectively the Levi-Civita connection, Riemann-Christoffel curvature tensor, Ricci tensor and scalar curvature of $M$. The curvature of a manifold plays the crucial role to determine the shape of the manifold. From a given metric one can determine the curvature but the converse is very cumbersome. For the sake of construction of a curvature restricted geometric structure one should impose a restriction on the curvature tensor by means of covariant derivatives or otherwise. It is well known that covariant derivative is a generalization of partial derivative and higher order of covariant derivatives imposed on a curvature tensor give rise different kinds of curvature restricted geometric structures. For example, the local symmetry and semisymmetry were introduced by Cartan (\cite{Cart26}, \cite{Cart46}) which are respectively appears due to the covariant derivative of first and second order. Again the locally symmetric manifold generalized by Chaki \cite{Chak87} as pseudosymmetric manifold. And Tam$\acute{\mbox{a}}$ssy and Binh's \cite{TB89} weakly symmetric manifold is another generalization of Chaki pseudosymmetric manifold. Also the notion of recurrent manifold appeared in the literature as a generalization of local symmetry. It may be noted that the study of recurrent manifold was initiated by Ruse (\cite{Ruse46}, \cite{Ruse49}, \cite{Ruse49a}) as kappa space and denoted as $K_n$. In 1979 Dubey \cite{Dube79} presented a generalization of $K_n$, called generalized recurrent manifold (briefly, $GK_n$) but we note that the structure $GK_n$ does not exist \cite{OO12} (see also \cite{Glod70}, \cite{Mike96}, \cite{MVH09}). Recently Shaikh and his coauthors introduced three notions of generalization of recurrent manifold, namely, quasi generalized recurrent manifold \cite{SR10} (briefly, $QGK_n$), hyper-generalized recurrent manifold \cite{SP10} (briefly, $HGK_n$) and weakly generalized recurrent manifold \cite{SR11} (briefly, $WGK_n$) along with their proper existence by suitable examples (see, \cite{SAR13}, \cite{SRK15}). We also note that for $\alpha =\beta$, a quasi-Einstein manifold ($S = \alpha g + \beta \eta\otimes\eta$) is $WGK_n$ if and only if it is $QGK_n$ and for $2\alpha =\beta$, a quasi-Einstein manifold is $HGK_n$ if and only if it is $QGK_n$. Very recently another generalization of such notion was given in \cite{SK14} and introduced the concept of super generalized recurrent manifold (briefly, $SGK_n$).\\
\indent The object of the present paper is to study the geometric properties of a $SGK_n$. Section 2 deals with the rudimentary facts of various curvature restricted geometric structures and tensors as preliminaries. Section 3 is devoted to the study of $SGK_n$ and contains the main results. It is proved that on a proper Roter type manifold \cite{Desz03} the notions of Ricci generalized recurrency \cite{DGK95} and super generalized recurrency are equivalent, and a $SGK_n$ satisfies semisymmetry condition if all of its associated 1-forms are closed and pairwise codirectional. The last section is concerned with the proper existence of such notion by means of a suitable example.
%
\section{Preliminaries}
Let us now consider a connected semi-Riemannian manifold $(M^n,g)$, $n\ge 3$. Let $C^{\infty}(M)$, $\chi(M)$, $\chi^*(M)$ and $\mathcal T^r_k(M)$ be respectively the algebra of all smooth functions, the Lie algebra of all smooth vector fields, the Lie algebra of all smooth 1-forms and the space of all smooth tensor fields of type $(r,k)$ on $M$. We now define some necessary terms and various curvature restricted geometric structures on $M$.\\
For $\Pi, \Phi \in \chi^*(M)$, we can define their exterior product $\Pi \wedge \Phi$ as:
$$\Pi \wedge \Phi = \frac{1}{2}\left(\Pi \otimes \Phi - \Phi \otimes \Pi\right),$$
where $\otimes$ denotes the tensor product. We note that if $\Pi \wedge \Phi = 0$, then $\Pi$ and $\Phi$ are linearly dependent or said to be codirectional. Again since $\nabla$ is torsion free so the exterior derivative $d\Pi$ of $\Pi$ can be expressed as:
$$d\Pi (X,Y) = (\nabla_X\Pi)(Y)-(\nabla_Y\Pi)(X)$$
for all $X,Y\in \chi(M)$. We also note that $\Pi$ is closed if $d\Pi =0$.\\
Now for $A,E\in \mathcal T^0_2(M)$ we have their Nomizu-Kulkarni product $A\wedge E$ as
\bea\label{eq2.1}
(A \wedge E)(X_1,X_2,Y_1,Y_2)&=&A(X_1,Y_2)E(X_2,Y_1) + A(X_2,Y_1)E(X_1,Y_2)\\\nonumber
&&-A(X_1,Y_1)E(X_2,Y_2) - A(X_2,Y_2)E(X_1,Y_1),
\eea
where $X_1, X_2, Y_1, Y_2\in \chi(M)$. Throughout the paper we consider $X, Y, X_i, Y_i \in \chi(M)$, $i = 1,2,\ldots$. As there is no confusion, here we use the same symbol $\wedge$ for Nomizu-Kulkarni product and exterior product.\\
Again for a symmetric $(0,2)$-tensor $A$ and $X, Y\in\chi(M)$ we can define two endomorphisms $\mathscr A$ and $X\wedge_{A} Y$ on $\chi(M)$ as:
$$g(\mathscr A X,Y) = A(X,Y) \ \ \mbox{and }\ (X \wedge_A Y)X_1 = A(Y,X_1)X - A(X,X_1)Y.$$
Then we get another $(0,2)$-tensor $A^2$, called the second level of $A$ with corresponding endomorphism operator $\mathscr A^2$ given below:
$$A^2(X,Y) = A(\mathscr A X,Y) = g(\mathscr A^2 X,Y).$$
We note that the endomorphisms $\mathscr A$, $\mathscr A^2$ and $X\wedge_A Y$are all $C^{\infty}(M)$-linear.
In particular we get the second level Ricci tensor $S^2$ given by
$$S^2(X,Y) = S(\mathscr S X,Y),$$
where $\mathscr S$ is the Ricci operator defined as $S(X,Y) = g(\mathscr SX,Y)$.\\
Using this Nomizu-Kurkarni product and $\wedge_{A}$ we can define some useful curvature tensors, namely, conformal curvature tensor $C$, projective curvature tensor $P$, concircular curvature tensor $W$ and conharmonic curvature tensor $K$ as follows:
$$C= R -\frac{1}{n-2} g \wedge S+\frac{\kappa}{2(n-1)(n-2)}g \wedge g,$$
$$P= R -\frac{1}{n-1} (\wedge_S),$$
$$W = R - \frac{\kappa}{2 n(n-1)}g \wedge g,$$
$$K = R - \frac{1}{n-2} g \wedge S.$$
Now for $D\in \mathcal T^0_4(M)$ we get its corresponding $\mathcal D\in \mathcal T^1_3(M)$ and the $C^{\infty}(M)$-linear endomorphism operator $\mathscr D(X_1,X_2)$ due to two vector fields $X$ and $Y$, as follows:
$$g(\mathcal D(X_1,X_2)X_3,X_4) = D(X_1,X_2,X_3,X_4) \ \mbox{ and}$$
$$\mathscr D(X,Y)X_3 = \mathcal D(X,Y)X_3.$$
We note that one can easily operate a $C^{\infty}(M)$-linear endomorphism $\mathscr L$ on $T\in\mathcal T^0_k(M)$ as follows:
\beb\label{rdot}
(\mathscr{L} T)(X_1,X_2,\cdots,X_k) = -T(\mathscr{L}X_1,X_2,\cdots,X_k) - \cdots - T(X_1,X_2,\cdots,\mathscr{L}X_k).
\eeb
Then for the endomorphisms $\mathscr D(X_1,X_2)$ and $X\wedge_{A} Y$ we get two $(0,k+2)$-tensors for $T$ as follows:
\beb\label{ddt}
&&(\mathscr D(X,Y) T)(X_1,X_2,\cdots,X_k)\\\nonumber
&&=-T(\mathscr D(X,Y)(X_1),X_2,\cdots,X_k) - \cdots - T(X_1,X_2,\cdots,\mathscr D(X,Y)(X_k))\\\nonumber
&&=- T(\mathcal D(X,Y)X_1,X_2,\cdots,X_k) - \cdots - T(X_1,X_2,\cdots, \mathcal D(X,Y)X_k),
\eeb
\beb\label{qgr}
&&((X \wedge_A Y) T)(X_1,X_2,\cdots,X_k)\\\nonumber
&&=-T((X \wedge_A Y)X_1,X_2,\cdots,X_k) - \cdots - T(X_1,X_2,\cdots,(X \wedge_A Y)X_k)\\\nonumber
&&= A(X, X_1) T(Y,X_2,\cdots,X_k) + \cdots + A(X, X_k) T(X_1,X_2,\cdots,Y)\\\nonumber
&&- A(Y, X_1) T(X,X_2,\cdots,X_k) - \cdots - A(Y, X_k) T(X_1,X_2,\cdots,X).
\eeb
We note that throughout this paper we denote\\ $(\mathscr D(X,Y) T)(X_1,X_2,\cdots,X_k)$ as $D\cdot T(X_1,X_2,\cdots,X_k,X,Y)$ and\\ $((X \wedge_A Y) T)(X_1,X_2,\cdots,X_k)$ as $Q(A,T)(X_1,X_2,\cdots,X_k,X,Y)$.\\
Again generalizing the notions of $A\wedge E$ and $\wedge_A$ for $A, E\in \mathcal T_2^0(M)$, we have the following for higher order tensors:
\beb
(A \wedge T)(X_1,X_2,Y_1,Y_2,\cdots,Y_k)&=&A(X_1,Y_2)T(X_2,Y_1,\cdots,Y_k) + A(X_2,Y_1)T(X_1,Y_2,\cdots,Y_k)\\
&-&A(X_1,Y_1)T(X_2,Y_2,\cdots,Y_k) - A(X_2,Y_2)T(X_1,Y_1,\cdots,Y_k).
\eeb
$$(X\wedge_T Y)(X_1,X_2,\cdots,X_k) = T(Y,X_1,X_3,\cdots X_{k})g(X,X_2) - T(X,X_1,X_3,\cdots X_{k})g(Y,X_2),$$
where $A\in \mathcal T_2^0(M)$ and $T\in \mathcal T_k^0(M)$.
Now from the definitions we can state the following:
\begin{pr}\label{pr2.1}
For $A\in\mathcal T^0_2(M)$ and $D\in\mathcal T^0_4(M)$, the following relations hold:\\
(i) $\nabla (X\wedge_A Y) = X\wedge_{\nabla A} Y$ and $\nabla(g\wedge A) = g\wedge (\nabla A)$,\\
(ii) $D\cdot (X\wedge_A Y) = X\wedge_{D\cdot A} Y$ and $D\cdot (g\wedge A) = g\wedge (D\cdot A)$ if $D\cdot g =0$.
\end{pr}
We now define some basic curvature restricted geometric structures:
\begin{defi}\label{rec}
Let $T\in \mathcal T^0_k(M)$. Then $M$ is said to be $T$-recurrent \cite{Ruse49} if the condition 
\be\label{kn}
\nabla T = \Pi \otimes T
\ee
holds on $\{x \in M:\nabla T\neq 0 \mbox{ at $x$}\}  \subset M$ for an 1-form $\Pi$, called the associated 1-form of this structure. Such an $n$-dimensional manifold is denoted by $T$-$K_{n}$ with associated 1-form $\Pi$ or simply $T$-$K_{n}$ with $\Pi$ or more simply $T$-$K_{n}$.
\end{defi}
\begin{defi}
For $T\in \mathcal T^0_4(M)$, $M$ is said to be $T$-quasi generalized recurrent \cite{SR10}, $T$-hyper generalized recurrent \cite{SP10} and $T$-weakly generalized recurrent \cite{SR11} respectively if the condition
\be\label{qgkn}
\nabla T = \Pi \otimes T + \Psi \otimes \left[g\wedge (g + \eta\otimes\eta)\right],
\ee
\be\label{hgkn}
\nabla T = \Pi \otimes T + \Psi \otimes g\wedge S \ \ \mbox{and}
\ee
\be\label{wgkn}
\nabla T = \Pi \otimes T + \Psi \otimes S\wedge S
\ee
holds respectively on $\{x\in M:\nabla T \neq \xi \otimes T \,\,\mbox{at}\,\, x\ \forall\ \xi \in \chi^*(M)\}  \subset M$ for some $\Pi$, $\Psi$ and $\eta \in \chi^*(M)$, called the associated 1-forms. An $n$-dimensional manifold satisfying \eqref{qgkn} $($resp., \eqref{hgkn} and \eqref{wgkn}$)$ is denoted by $T$-$QGK_{n}$ $($resp., $T$-$HGK_{n}$ and $T$-$WGK_{n})$.
\end{defi}
We note that we call these structures in short as recurrent like structures. Now Generalizing these recurrent like structures $T$-$WGK_n$ and $T$-$HGK_n$ we now define the super generalized recurrent structure on $M$.
\begin{defi}
For $T\in \mathcal T^0_4$, $M$ is said to be $T$-super generalized recurrent manifold \cite{SK14} if the following condition
\be\label{tsgk}
\nabla T = \Pi \otimes T + \Phi \otimes S\wedge S + \Psi \otimes g\wedge S + \Theta \otimes g\wedge g
\ee
holds on $\{x\in M:\nabla T - \xi \otimes T - \zeta \otimes S\wedge S - \theta \otimes g\wedge S \neq 0 \,\,\mbox{at}\,\, x\ \forall\ \xi, \zeta, \theta \in \chi^*(M)\}  \subset M$ for some 1-forms $\Pi$, $\Phi$, $\Psi$ and $\Theta$, called the associated 1-forms of this structure. Such an $n$-dimensional manifold is denoted by $T$-$SGK_{n}$ with associated 1-forms $(\Pi, \Phi, \Psi, \Theta)$ or simply $T$-$SGK_{n}$ with $(\Pi, \Phi, \Psi, \Theta)$ or more simply $T$-$SGK_n$.
\end{defi}
We note that if we take the particular value of $T$ as the Riemann-Christoffel curvature tensor $R$, then we call $R$-$K_n$ as simply recurrent manifold and denoted as $K_n$. Similarly we call $R$-$QGK_n$, $R$-$HGK_{n}$, $R$-$WGK_{n}$ and $R$-$SGK_{n}$ as quasi generalized recurrent manifold, hyper generalized recurrent manifold, weakly generalized recurrent manifold and super generalized recurrent manifold respectively, and denoted them as simply $QGK_{n}$, $HGK_{n}$, $WGK_{n}$ and $SGK_{n}$ respectively. Again if we take $T$ as the Ricci tensor $S$, then we call $S$-$K_n$ as Ricci recurrent manifold.\\
\indent Again to express the defining condition of $T$-$QGK_n$ in more explicit way we get its extension as
$$\nabla T = \Pi \otimes T + \Phi \otimes g\wedge g + \Psi\otimes g\wedge(\eta\otimes\eta),$$
where $\Pi, \Phi, \Psi \in \chi^*(M)$, and such an $n$-dimensional manifold is called as quasi generalized recurrent like.\\
Again there is a generalization of the notion of recurrency for $(0,2)$-tensors as follows:
\begin{defi}
For $Z\in \mathcal T^0_2$, $M$ is said to be $Z$-generalized recurrent manifold \cite{DGK95} $($briefly, $Z$-$GK_n$ $)$ if the following condition
\be\label{tgk}
\nabla Z = \Pi \otimes Z + \Phi \otimes g
\ee
holds on $\{x\in M:\nabla Z \neq \xi \otimes Z \,\,\mbox{at}\,\, x\ \forall\ \xi \in \chi^*(M)\}  \subset M$ for some $\Pi$ and $\Phi\in\chi^*(M)$, called the associated 1-forms of this structure.
\end{defi}
In particular if $Z=S$, then we get Ricci generalized recurrent $($briefly, $S$-$GK_n$ $)$.
\begin{defi}
For $T\in\mathcal T^0_k$, $M$ is said to be $T$-semisymmetric $($briefly, $T$-$SS_n$ $)$ $($\cite{Cart46}, \cite{Szab82}$)$ if
$$R\cdot T = 0.$$
Again $M$ is said to be $T$-pseudosymmetric in sense of Deszcz $($briefly, $T$-$PS_n$ $)$ $($\cite{AD83}, \cite{Desz92}$)$ if $R\cdot T$ and $Q(g,T)$ are linearly dependent, i.e.,
$$R\cdot T = L_T Q(g,T)$$
holds on $\{x\in M: Q(g,T) \neq 0 \mbox{ at $x$}\}$ for some function $L_T$ on the set.
\end{defi}
In particular, $R$-$SS_n$ and $R$-$PS_n$ are respectively known as, simply, semisymmetric manifold $($briefly, $SS_n$ $)$ and pseudosymmetric manifold $($briefly, $PS_n$ $)$.
\begin{defi}
The manifold $M$ is said to be Roter type $($briefly, $RT_n$ $)$ $($\cite{Desz03}, \cite{Desz03a}, \cite{DGHS11}$)$ if its curvature tensor $R$ has the following decomposition:
$$R = N_1 g\wedge g + N_2 g\wedge S + N_3 S\wedge S,$$
for some $N_1, N_2$ and $N_3\in C^{\infty}(M)$. Moreover it is said to be proper Roter type manifold if $N_3 \neq 0$.
\end{defi}
We note that recently Shaikh and his coauthors introduced a generalized notion of Roter type structure \cite{SDHJK15}. For more details about generalized Roter type manifold and its characterization on a warped product manifold we refer the readers to see \cite{SKgrt} and \cite{SKgrtw} and also references therein.\\
\indent Now contracting the decomposition relation of $R$ on a $RT_n$, we get a generalization of Einstein manifold, namely $Ein(2)$. $M$ is said to be Einstein (resp., $Ein(2)$) if
$$S=\frac{\kappa}{n}g \ \ \ \mbox{(resp., $a_1 S^2 + a_2 S + a_3 g = 0$),}$$
where $a_1, a_2, a_3\in C^{\infty}(M)$.
We note that there is an another generalization of Einstein manifold, namely, quasi-Einstein manifold. The manifold $M$ is said to be quasi-Einstein if
\be\label{qEn}
S = \alpha g +\beta \eta\otimes\eta
\ee
holds for some $\alpha, \beta \in C^{\infty}(M)$ and $\eta\in \chi^*(M)$.
%
\section{Main results}
Let us consider a super generalized recurrent condition on $M$ as
\be\label{sgk}
\nabla R = \Pi \otimes R + \Phi \otimes S\wedge S + \Psi \otimes g\wedge S + \Theta \otimes g\wedge g.
\ee
Then contracting above equation we get
\be\label{csgk}
\nabla S = \Pi_1\otimes S^2 + \Phi_1\otimes S + \Psi_1\otimes g,
\ee
where $\Pi_1 = -2\Phi$, $\Phi_1 = \Pi + 2\kappa \Phi + (n-2)\Psi$ and $\Psi_1 = \kappa\Psi + 2(n-1)\Theta$.
\begin{thm}\label{th3.1}
Let $M$ be a $SGK_n$ satisfying \eqref{sgk}. Then the associated 1-forms are not uniquely determined.
\end{thm}
\noindent \textbf{Proof:} We know that $R$ satisfies the second Bianchi identity
$$(\nabla_{X_1})R(X_2,X_3,X_4,X_5)+(\nabla_{X_2})R(X_3,X_1,X_4,X_5)+(\nabla_{X_3})R(X_1,X_2,X_4,X_5) = 0.$$
Then putting the value of $\nabla R$ from \eqref{sgk}, we get
\bea\label{eq3.2}
&&\sum_{X_1,X_2,X_3}\Big[\Pi(X_1)R(X_2,X_3,X_4,X_5)+\Phi(X_1)(S\wedge S)(X_2,X_3,X_4,X_5)\\\nonumber
&&\hspace{0.5in}+\Psi(X_1)(g\wedge S)(X_2,X_3,X_4,X_5)+\Theta(X_1)(g\wedge g)(X_2,X_3,X_4,X_5)\Big]=0,
\eea
where $\sum\limits_{\scriptscriptstyle{X_1,X_2,X_3}}$ denotes the cyclic sum in $X_1,X_2$ and $X_3$. Now contracting above in $X_1$ and $X_5$, we get
\beb
&&-R(V, X_4, X_2, X_3)+\{\kappa \Psi(X_3) - \Psi(\mathcal S(X_3)) + 2 (n-2) \Theta(X_3)\}g(X_2, X_4)\\
&&+\{-\kappa \Psi(X_2) + \Psi(\mathcal S(X_2)) - 2 (n-2) \Theta(X_2)\}g(X_3, X_4)\\
&&+\{\Pi(X_3) + 2 \kappa \Phi(X_3) - 2 \Phi(\mathcal S(X_3)) + (n-3) \Psi(X_3)\}S(X_2, X_4)-2 \Phi(X_3)S^2(X_2, X_4)\\
&&+\{-\Pi(X_2) - 2 \kappa \Phi(X_2) + 2 \Phi(\mathcal S(X_2)) - (n-3) \Psi(X_2)\}S(X_3, X_4)+2 \Phi(X_2)S^2(X_3, X_4) = 0,
\eeb
where $V$ is the corresponding vector field of $\Pi$, i.e., $g(V,X) = \Pi(X)$, for all $X\in\chi(M)$.
Again contracting in $X_3$ and $X_4$, we get
\beb
&&-\kappa \Pi(X_2) +  2 \Big[\Pi(\mathcal S(X_2)) + (\kappa^{(2)} - \kappa^2) \Phi(X_2) + 2 \kappa \Phi(\mathcal S(X_2))\\
&& - 2 \Phi(\mathcal S^2(X_2)) - (n-2) \{\kappa \Psi(X_2) - \Psi(\mathcal S(X_2)) + (n-1) \Theta(X_2)\}\Big]=0,
\eeb
where $\kappa^{(2)}$ is the trace of $S^2$. Hence the result.\\
\noindent \textbf{Note:} We note that the result is true for $C$-$SGK_n$, $P$-$SGK_n$, $W$-$SGK_n$ and $K$-$SGK_n$ also.
\begin{thm}
Let $M$ be a $SGK_n$ satisfying \eqref{sgk}. Then the associated 1-forms are linearly dependent with $d\kappa$ such that
$$d\kappa = \kappa \Pi + 2(\kappa^2 - \kappa^{(2)}) \Phi + 2(n-1)[\kappa \Psi + n \Theta],$$
where $\kappa^{(2)} = Tr(S^2)$.
\end{thm}
\noindent \textbf{Proof:} Since the manifold satisfies \eqref{sgk} so it satisfies \eqref{csgk}, i.e.,
\bea
\nabla_X S(X_1,X_2) = -2\Phi(X) S^2(X_1,X_2) &+& (\Pi + 2\kappa \Phi + (n-2)\Psi)(X) S(X_1,X_2)\\\nonumber
&+& (\kappa\Psi + 2(n-1)\Theta)(X) g(X_1,X_2).
\eea
Now contracting the above equation in $X_1$ and $X_2$, we get the results easily.
\begin{thm}\label{en}
A super generalized recurrent manifold becomes a recurrent manifold if it is Einstein and in this case the relation between the associated 1-forms are given by
$$\frac{\kappa^2}{n^2}\Phi + \frac{\kappa}{n}\Psi + \Theta = 0.$$
\end{thm}
\noindent \textbf{Proof:} Since the manifold is Einstein, so $S = \frac{\kappa}{n}g$ and thus from \eqref{sgk} we have:
$$\nabla R = \Pi \otimes R + \left[\frac{\kappa^2}{n^2}\Phi + \frac{\kappa}{n}\Psi + \Theta\right] \otimes g\wedge g$$
Again in \cite{OO12} Olszak and Olszak showed that for any semi-Riemannian manifold satisfying such curvature condition, sometime called generalized recurrent manifold \cite{Dube79}, the coefficient of $g\wedge g$ is zero. Thus we can conclude that the manifold becomes recurrent and $\frac{\kappa^2}{n^2}\Phi + \frac{\kappa}{n}\Psi + \Theta = 0$.
\begin{thm}
Let $M$ be a $SGK_n$ with $(\Pi, \Phi, \Psi, \Theta)$ and also quasi-Einstein satisfying \eqref{qEn}, then it becomes a $QGK_n$-like with $(\Pi,\alpha^2\Phi+\alpha\Psi+\Theta,2\alpha\beta\Phi + \beta\Psi, \eta)$. Moreover it becomes a $QGK_n$ if
$$\alpha^2\Phi+\alpha\Psi+\Theta= \gamma(2\alpha\beta\Phi + \beta\Psi),$$
where $\gamma \in C^{\infty}(M)$ and $\gamma> 0$ on $M$.
\end{thm}
\noindent \textbf{Proof:} The proof is similar to the proof of the Theorem \ref{en}.\\
Now we evaluate the conditions on the associated 1-forms such that the $SGK_n$ satisfies semisymmetric condition.
\begin{thm}
Let $M$ be a $SGK_n$. Then it satisfies semisymmetry condition if all $\Pi, \Phi, \Psi$ and $\Theta$ are closed and pairwise codirectional.
\end{thm}
\noindent \textbf{Proof:} Let $M$ be a $SGK_n$ with $(\Pi, \Phi, \Psi, \Theta)$. Then
$$\nabla_X R = \Pi(X) \otimes R + \Phi(X) \otimes S\wedge S + \Psi(X) \otimes g\wedge S + \Theta(X) \otimes g\wedge g$$
Again differentiating covariantly with respect to $Y$, we get
\beb
&&\nabla_Y(\nabla_X R) = (\nabla_Y \Pi)(X) R + (\nabla_Y \Phi)(X) S\wedge S + (\nabla_Y \Psi)(X) g\wedge S + (\nabla_Y \Theta)(X) g\wedge g\\
&&+\left[\Pi(X) \Psi(Y)+\Psi(X) \Phi _1(Y)+2 \Phi(X) \Psi _1(Y)\right] g\wedge S\\
&&+ \Pi(X)\Pi(Y) R + \Psi(X) \Pi _1(Y) g\wedge S^2 + \left[\Pi(X) \Theta(Y)+\Psi(X) \Psi _1(Y)\right] g\wedge g\\
&&+2 \Phi(X) \Pi _1(Y) S\wedge S^2 + \left[\Pi(X) \Phi(Y)+2 \Phi(X) \Phi _1(Y)\right] S\wedge S.
\eeb
Now from definition we have
\beb
&& R(X,Y)\cdot R = \nabla_X(\nabla_Y R) - \nabla_Y(\nabla_X R)\\
&&= d\Pi(X,Y) R + 2\left[\Phi(Y) \Psi(X) - \Phi(X) \Psi(Y)\right] g\wedge S^2\\
&&+ \left[d\Phi(X,Y)+\Phi(Y)(\Pi(X)+2 (n-2) \Psi(X))-\Phi(X) (\Pi(Y)+2 (n-2) \Psi(Y))\right] S\wedge S \\
&&+ \left[d\Psi(X,Y) - 4 (n-1)\left(\Theta(Y) \Phi(X) - \Theta(X) \Phi(Y)\right)\right] g\wedge S \\
&&+ \left[d\Theta(X,Y)-\Theta(Y) (\Pi(X)+2 (n-1) \Psi(X))+\Theta(X) (\Pi(Y)+2 (n-1) \Psi(Y))\right] g\wedge g.
\eeb
Thus we have
\bea\label{eq3.5}
R\cdot R &=& d\Pi R - 4 (\Phi\wedge\Psi) g\wedge S^2 + [d\Phi - 2\Phi\wedge(\Pi - 2(n-2)\Psi)] S\wedge S\\\nonumber
&+& [d\Psi - 8(n-1)\Phi\wedge\Theta] g\wedge S + [d\Theta - 2\Theta\wedge(\Pi+2(n-1)\Psi)] g\wedge g.
\eea
Thus from above equation we can easily conclude the required result.
\begin{pr}\label{pr3.1}
A super generalized recurrent manifold $M$ satisfying \eqref{sgk} is Ricci generalized recurrent if and only if $\Phi =0$ or $M$ satisfies some proper $Ein(2)$ condition.
\end{pr}
\noindent \textbf{Proof:} The result easily follows from \eqref{csgk}.
\begin{thm}\label{th3.7}
On a proper Roter type manifold the notion of Ricci generalized recurrency and super generalized recurrency are equivalent.
\end{thm}
\noindent \textbf{Proof:} Let us consider the proper Roter type condition
$$R = N_1 g\wedge g + N_2 g\wedge S + N_3 S\wedge S,$$
where $N_1, N_2$ and $N_3 (\ne 0)\in\chi^*(M)$. Then obviously the manifold satisfies some $Ein(2)$ condition and thus by Proposition \ref{pr3.1}, we conclude that the manifold is Ricci generalized recurrent if it is super generalized recurrent.
Again differentiating covariantly the Roter type condition, we get the inverse implication.
\begin{thm}\label{th3.8}
Let $M$ be a $S$-$GK_n$ with $(\overline\Pi, \overline\Phi)$ and also $SGK_n$ with $(\Pi, \Phi, \Psi, \Theta)$. Then $M$ is\\
(i) $C$-$SGK_n$ with $(\Pi, \Phi, \Psi-\frac{\overline\Pi-\Pi}{n-2}, \Theta-\frac{\overline\Phi}{n-2}+\frac{\kappa \overline\Pi + n \overline\Phi-\kappa\Pi}{2(n-1)(n-2)})$.\\
(ii) $W$-$SGK_n$ with $(\Pi, \Phi, \Psi, \Theta-\frac{\kappa \overline\Pi + n \overline\Phi-\kappa\Pi}{2n(n-1)})$.\\
(iii) $K$-$SGK_n$ with $(\Pi, \Phi, \Psi-\frac{\overline\Pi-\Pi}{n-2}, \Theta-\frac{\overline\Phi}{n-2})$. And, especially,\\
(iv) $P$-$SGK_n$ with $(\Pi, \Phi, \Psi, \Theta-\frac{\overline\Phi}{2(n-2)})$ if $\Pi = \overline\Pi$.
\end{thm}
\noindent \textbf{Proof:} The results immediately follow from definition with help of Proposition \ref{pr2.1}.\\
\indent Again to get the condition for a $SGK_n$ to be $C$-$SGK_n$, $P$-$SGK_n$, $W$-$SGK_n$ or $K$-$SGK_n$ with same associated 1-forms, we have the following:
\begin{thm}\label{th3.9}
Let $T\in \mathcal T^0_4(M)$. Then the following conditions are equivalent:\\
(i) the tensor field $(T-R)\in \mathcal T^0_4(M)$ is recurrent with $\Pi$.\\
(ii) the structures $T$-$SGK_n$ with $(\Pi, \Phi, \Psi, \Theta)$ and $SGK_n$ with $(\Pi, \Phi, \Psi, \Theta)$ are equivalent.
\end{thm}
\noindent \textbf{Proof:} The result easily follows from definitions.\\
\indent From above we have the following interesting results.
\begin{cor}\label{cor3.1}
Let $M$ be a $SGK_n$ with $(\Pi, \Phi, \Psi, \Theta)$. Then it is\\
(i) $C$-$SGK_n$ with $(\Pi, \Phi, \Psi, \Theta)$ if and only if it is Ricci recurrent with $\Pi$.\\
(ii) $P$-$SGK_n$ with $(\Pi, \Phi, \Psi, \Theta)$ if and only if it is Ricci recurrent with $\Pi$.\\
(iii) $W$-$SGK_n$ with $(\Pi, \Phi, \Psi, \Theta)$ if and only if $d\kappa = \kappa\Pi$.\\
(iv) $K$-$SGK_n$ with $(\Pi, \Phi, \Psi, \Theta)$ if and only if it is Ricci recurrent with $\Pi$.
\end{cor}
\noindent \textbf{Proof:} (i) From Theorem \ref{th3.9}, $M$ is $C$-$SGK_n$ if and only if the $\frac{1}{n-2} g \wedge S-\frac{\kappa}{2(n-1)(n-2)}g \wedge g = \frac{1}{n-2}\left[g\wedge (S-\frac{2\kappa}{n-1}g)\right]$ is recurrent with $\Pi$, or equivalently the tensor $(S-\frac{2\kappa}{n-1}g)$ is recurrent with $\Pi$ (by Proposition \ref{pr2.1}), i.e.,
$$\nabla S - 2 \frac{d\kappa}{n-1} g = \Pi\otimes\left(S - 2 \frac{\kappa}{n-1} g\right).$$
Now the above condition implies that $d\kappa = \kappa\Pi$ and thus Ricci recurrent with $\Pi$. Again Ricci recurrent with $\Pi$ implies the above condition. Hence (i) is proved.\\
(ii) Also from Theorem \ref{th3.9}, $M$ is $P$-$SGK_n$ if and only if $\frac{1}{n-2} X \wedge_S Y$ is recurrent, or equivalently $S$ is recurrent (By Proposition \ref{pr2.1}).\\
The proof of (iii) and (iv) are similar to the above.\\
\noindent \textbf{Note:} The sufficient parts of the above corollary can be directly reduced from Theorem \ref{th3.8} by using $\overline\Pi =\Pi$ and $\overline\Phi =0$.
\section{\bf Examples}\label{exam}
\textbf{Example 1:} Consider the open connected subset $M= \{(x^1,x^2,x^3,x^4)\in\mathbb R^4:x^1,x^2,x^3,x^4> 0\}$ of $\mathbb R^4$ such that on $M$
$$(x^1)^4 (x^3)^2 (x^2)^4+(x^1)^4 (x^4)^2 (x^2)^4+2 (x^1)^4 x^3 x^4 (x^2)^4-\left(x^1+x^2\right)^2 (x^3)^4 (x^4)^4 > 0.$$
Take the following Riemannian metric on $M$:
$$ds^2 = x^2(dx^1)^2+x^1(dx^2)^2+x^4(dx^3)^2+x^3(dx^4)^2.$$
We can easily evaluate the local components of various tensors on $M$. The non-zero components (upto symmetry) of the Riemann-Christoffel curvature tensor $R$ and Ricci tensor $S$ are given by
$$ R_{1212}=\frac{1}{4} \left(\frac{1}{x^2}+\frac{1}{x^1}\right), \ R_{3434}=\frac{1}{4} \left(\frac{1}{x^4}+\frac{1}{x^3}\right) \mbox{ and}$$
$$S_{11}=-\frac{\frac{x^1}{x^2}+1}{4 (x^1)^2}, \ S_{22}=-\frac{\frac{x^2}{x^1}+1}{4 (x^2)^2}, \ S_{33}=-\frac{\frac{x^3}{x^4}+1}{4 (x^3)^2}, \ S_{44}=-\frac{\frac{x^4}{x^3}+1}{4 (x^4)^2}.$$
Again the non-zero local components (upto symmetry) of $\nabla R$ are
$$ R_{1212,1}=-\frac{\frac{x^1}{x^2}+2}{4 (x^1)^2}, \  R_{1212,2}=-\frac{\frac{x^2}{x^1}+2}{4 (x^2)^2}, \ R_{3434,3}=-\frac{\frac{x^3}{x^4}+2}{4 (x^3)^2}, \  R_{3434,4}=-\frac{\frac{x^4}{x^3}+2}{4 (x^4)^2}.$$
Now we have the non-zero local components (upto symmetry) of $g\wedge g$, $g\wedge S$ and $S\wedge S$ as follows:
$$ g\wedge g_{1212}=-2x^1 x^2, \ g\wedge g_{1313}=-2x^2 x^4, \ g\wedge g_{1414}=-2x^2 x^3,$$
$$g\wedge g_{2323}=-2x^1 x^4, \ g\wedge g_{2424}=-2x^1 x^3, \ g\wedge g_{3434}=-2x^3 x^4;$$
$$ g\wedge S_{1212}=\frac{1}{2} \left(\frac{1}{x^2}+\frac{1}{x^1}\right), \ g\wedge S_{1313}=\frac{\left(x^1+x^2\right) x^4}{4 (x^1)^2 x^2}+\frac{x^2 \left(x^3+x^4\right)}{4 (x^3)^2 x^4},$$
$$ g\wedge S_{1414}=\frac{\left(x^1+x^2\right) x^3}{4 (x^1)^2 x^2}+\frac{x^2 \left(x^3+x^4\right)}{4 (x^4)^2 x^3}, \ g\wedge S_{2323}=\frac{\left(x^1+x^2\right) x^4}{4 x^1 (x^2)^2}+\frac{x^1 \left(x^3+x^4\right)}{4 (x^3)^2 x^4},$$
$$ g\wedge S_{2424}=\frac{\left(x^1+x^2\right) x^3}{4 x^1 (x^2)^2}+\frac{x^1 \left(x^3+x^4\right)}{4 (x^4)^2 x^3}, \ g\wedge S_{3434}=\frac{1}{2} \left(\frac{1}{x^4}+\frac{1}{x^3}\right) \mbox{ and}$$
$$ S\wedge S_{1212}=-\frac{\left(x^1+x^2\right)^2}{8 (x^1)^3 (x^2)^3}, \ S\wedge S_{1313}=-\frac{\left(x^1+x^2\right) \left(x^3+x^4\right)}{8 (x^1)^2 x^2 (x^3)^2 x^4},$$
$$ S\wedge S_{1414}=-\frac{\left(x^1+x^2\right) \left(x^3+x^4\right)}{8 (x^1)^2 x^2 x^3 (x^4)^2}, \ S\wedge S_{2323}=-\frac{\left(x^1+x^2\right) \left(x^3+x^4\right)}{8 x^1 (x^2)^2 (x^3)^2 x^4},$$
$$ S\wedge S_{2424}=-\frac{\left(x^1+x^2\right) \left(x^3+x^4\right)}{8 x^1 (x^2)^2 x^3 (x^4)^2}, \ S\wedge S_{3434}=-\frac{\left(x^3+x^4\right)^2}{8 (x^3)^3 (x^4)^3}.$$
We can now test that the manifold $M$ satisfies the super generalized recurrent condition as:
$$\nabla R = \Pi \otimes R + \Phi \otimes (S\wedge S) + \Psi \otimes (g\wedge S) + \Theta \otimes g\wedge g$$
where $\Pi$, $\Phi$, $\Psi$ and $\Theta$ are given by
\be
\Pi_{i}=\left\{\begin{array}{ccc}
&\frac{8 \Theta_1 \left((x^1)^2 x^3 (x^2)^2+(x^1)^2 x^4 (x^2)^2-\left(x^1+x^2\right) (x^3)^2 (x^4)^2\right)^2-x^1 (x^2)^2 \left(x^1+2 x^2\right) \left(x^3+x^4\right)^2}{\left(x^1+x^2\right) \left(x^3+x^4\right) \left((x^1)^2 \left(x^3+x^4\right) (x^2)^2+(x^3)^2 (x^4)^2 x^2+x^1 (x^3)^2 (x^4)^2\right)}&\mbox{for} \ i = 1\\

&\frac{8 \Theta_2 \left((x^1)^2 x^3 (x^2)^2+(x^1)^2 x^4 (x^2)^2-\left(x^1+x^2\right) (x^3)^2 (x^4)^2\right)^2-(x^1)^2 x^2 \left(2 x^1+x^2\right)\left(x^3+x^4\right)^2}{\left(x^1+x^2\right) \left(x^3+x^4\right) \left((x^1)^2 \left(x^3+x^4\right) (x^2)^2+(x^3)^2 (x^4)^2 x^2+x^1 (x^3)^2 (x^4)^2\right)}&\mbox{for} \ i = 2\\

&\frac{8 \Theta_3 \left((x^1)^2 x^3 (x^2)^2+(x^1)^2 x^4 (x^2)^2-\left(x^1+x^2\right) (x^3)^2 (x^4)^2\right)^2-\left(x^1+x^2\right)^2 x^3 (x^4)^2 \left(x^3+2 x^4\right)}{\left(x^1+x^2\right) \left(x^3+x^4\right) \left((x^1)^2 \left(x^3+x^4\right) (x^2)^2+(x^3)^2 (x^4)^2 x^2+x^1 (x^3)^2 (x^4)^2\right)}&\mbox{for} \ i = 3\\

&\frac{8 \Theta_4 \left((x^1)^2 x^3 (x^2)^2+(x^1)^2 x^4 (x^2)^2-\left(x^1+x^2\right) (x^3)^2 (x^4)^2\right)^2-\left(x^1+x^2\right)^2 (x^3)^2 x^4 \left(2 x^3+x^4\right)}{\left(x^1+x^2\right) \left(x^3+x^4\right) \left((x^1)^2 \left(x^3+x^4\right) (x^2)^2+(x^3)^2 (x^4)^2 x^2+x^1 (x^3)^2 (x^4)^2\right)}&\mbox{for} \ i = 4,\\
\end{array}\right.
\ee
\be
\Phi_{i}=\left\{\begin{array}{ccc}
&\frac{2 x^1 (x^2)^2 (x^3)^2 (x^4)^2 \left(\frac{8 \Theta_1 x^1}{x^3+x^4}-\frac{x^1+2 x^2}{(x^1)^2 x^3 (x^2)^2+(x^1)^2 x^4 (x^2)^2-\left(x^1+x^2\right) (x^3)^2 (x^4)^2}\right)}{x^1+x^2}&\mbox{for} \ i = 1\\

&\frac{2 (x^1)^2 x^2 (x^3)^2 (x^4)^2 \left(\frac{8 \Theta_2 x^2}{x^3+x^4}-\frac{2 x^1+x^2}{(x^1)^2 x^3 (x^2)^2+(x^1)^2 x^4 (x^2)^2-\left(x^1+x^2\right) (x^3)^2 (x^4)^2}\right)}{x^1+x^2}&\mbox{for} \ i = 2\\

&\frac{2 (x^1)^2 (x^2)^2 x^3 (x^4)^2 \left(\frac{8 \Theta_3 x^3}{x^1+x^2}+\frac{x^3+2 x^4}{(x^1)^2 x^3 (x^2)^2+(x^1)^2 x^4 (x^2)^2-\left(x^1+x^2\right) (x^3)^2 (x^4)^2}\right)}{x^3+x^4}&\mbox{for} \ i = 3\\

&\frac{2 (x^1)^2 (x^2)^2 (x^3)^2 x^4 \left(\frac{8 \Theta_4 x^4}{x^1+x^2}+\frac{2 x^3+x^4}{(x^1)^2 x^3 (x^2)^2+(x^1)^2 x^4 (x^2)^2-\left(x^1+x^2\right) (x^3)^2 (x^4)^2}\right)}{x^3+x^4}&\mbox{for} \ i = 4,\\
\end{array}\right.
\ee
\be
\Psi_{i}=\left\{\begin{array}{ccc}
&\frac{x^1 (x^2 x^3 x^4)^2 \left[16 \Theta_1 x^1 \left((x^1)^2 x^3 (x^2)^2+(x^1)^2 x^4 (x^2)^2-\left(x^1+x^2\right) (x^3)^2 (x^4)^2\right)-\left(x^1+2 x^2\right) \left(x^3+x^4\right)\right]}{(x^1)^4 (x^3)^2 (x^2)^4+(x^1)^4 (x^4)^2 (x^2)^4+2 (x^1)^4 x^3 x^4 (x^2)^4-\left(x^1+x^2\right)^2 (x^3)^4 (x^4)^4}&\mbox{for} \ i = 1\\

&\frac{x^2 (x^1 x^3 x^4)^2 \left[16 \Theta_2 x^2 \left((x^1)^2 x^3 (x^2)^2+(x^1)^2 x^4 (x^2)^2-\left(x^1+x^2\right) (x^3)^2 (x^4)^2\right)-\left(2 x^1+x^2\right) \left(x^3+x^4\right)\right]}{(x^1)^4 (x^3)^2 (x^2)^4+(x^1)^4 (x^4)^2 (x^2)^4+2 (x^1)^4 x^3 x^4 (x^2)^4-\left(x^1+x^2\right)^2 (x^3)^4 (x^4)^4}&\mbox{for} \ i = 2\\

&\frac{x^3 (x^1 x^2 x^4)^2 \left[16 \Theta_3 x^3 \left((x^1)^2 x^3 (x^2)^2+(x^1)^2 x^4 (x^2)^2-\left(x^1+x^2\right) (x^3)^2 (x^4)^2\right)+\left(x^1+x^2\right) \left(x^3+2 x^4\right)\right]}{(x^1)^4 (x^3)^2 (x^2)^4+(x^1)^4 (x^4)^2 (x^2)^4+2 (x^1)^4 x^3 x^4 (x^2)^4-\left(x^1+x^2\right)^2 (x^3)^4 (x^4)^4}&\mbox{for} \ i = 3\\

&\frac{x^4 (x^1 x^2 x^3)^2 \left[16 \Theta_4 x^4 \left((x^1)^2 x^3 (x^2)^2+(x^1)^2 x^4 (x^2)^2-\left(x^1+x^2\right) (x^3)^2 (x^4)^2\right)+\left(x^1+x^2\right) \left(2 x^3+x^4\right)\right]}{(x^1)^4 (x^3)^2 (x^2)^4+(x^1)^4 (x^4)^2 (x^2)^4+2 (x^1)^4 x^3 x^4 (x^2)^4-\left(x^1+x^2\right)^2 (x^3)^4 (x^4)^4}&\mbox{for} \ i = 4.\\
\end{array}\right.
\ee
We see that the associated $1$-forms are not uniquely determined as we can take $\Theta$ arbitrarily and others are related with it (which supports the Theorem \ref{th3.1}). We can check that the manifold is neither hyper generalized nor weakly generalized recurrent.\\
Again from the values of $R$, $g\wedge g$, $g\wedge S$ and $S\wedge S$, we can easily check that the manifold is Roter type satisfying $R = N_1 g\wedge g + N_2 g\wedge S + N_3 S\wedge S$, where
\beb
N_1 &=& -\frac{\left(x^1+x^2\right) \left(x^3+x^4\right) \left((x^1)^2 \left(x^3+x^4\right) (x^2)^2+(x^3)^2 (x^4)^2 x^2+x^1 (x^3)^2 (x^4)^2\right)}{8 \left(-(x^1)^2 (x^2)^2 \left(x^3+x^4\right) +(x^3)^2 (x^4)^2 (x^1 + x^2)\right)^2},\\
N_2 &=& - \frac{2 (x^1)^2 (x^2)^2 \left(x^1+x^2\right) (x^3)^2 (x^4)^2 \left(x^3+x^4\right)}{\left(-(x^1)^2 (x^2)^2 \left(x^3+x^4\right) +(x^3)^2 (x^4)^2 (x^1 + x^2)\right)^2} \ \mbox{ and}\\
N_3 &=& - \frac{2 (x^1)^2 (x^2)^2 (x^3)^2 (x^4)^2 \left((x^1)^2 \left(x^3+x^4\right) (x^2)^2+(x^3)^2 (x^4)^2 x^2+x^1 (x^3)^2 (x^4)^2\right)}{\left(-(x^1)^2 (x^2)^2 \left(x^3+x^4\right) +(x^3)^2 (x^4)^2 (x^1 + x^2)\right)^2}.\\
\eeb
Thus by Theorem \ref{th3.7}, the manifold is Ricci generalized recurrent satisfying:
$$\nabla S = \overline\Pi \otimes S + \overline\Phi \otimes g,$$
where $\overline\Pi$ and $\overline\Phi$ are given by
\be
\overline\Pi_{i}=\left\{\begin{array}{ccc}
&\frac{\left(x^1+2 x^2\right) (x^3)^2 (x^4)^2}{(x^2)^2 x^3 (x^1)^3+(x^2)^2 x^4 (x^1)^3-\left(x^1+x^2\right) (x^3)^2 (x^4)^2 x^1}&\mbox{for} \ i = 1\\

&\frac{\left(2 x^1+x^2\right) (x^3)^2 (x^4)^2}{(x^1)^2 x^3 (x^2)^3+(x^1)^2 x^4 (x^2)^3-\left(x^1+x^2\right) (x^3)^2 (x^4)^2 x^2}&\mbox{for} \ i = 2\\

&\frac{(x^1)^2 (x^2)^2 \left(x^3+2 x^4\right)}{x^3 \left(-(x^1)^2 \left(x^3+x^4\right) (x^2)^2+(x^3)^2 (x^4)^2 x^2+x^1 (x^3)^2 (x^4)^2\right)}&\mbox{for} \ i = 3\\

&\frac{(x^1)^2 (x^2)^2 \left(2 x^3+x^4\right)}{x^4 \left(-(x^1)^2 \left(x^3+x^4\right) (x^2)^2+(x^3)^2 (x^4)^2 x^2+x^1 (x^3)^2 (x^4)^2\right)}&\mbox{for} \ i = 4,\\
\end{array}\right.
\ee
\be
\overline\Phi_{i}=\left\{\begin{array}{ccc}
&\frac{\left(x^1+2 x^2\right) \left(x^3+x^4\right)}{4 (x^2)^2 x^3 (x^1)^3+4 (x^2)^2 x^4 (x^1)^3-4 \left(x^1+x^2\right) (x^3)^2 (x^4)^2 x^1}&\mbox{for} \ i = 1\\

&\frac{\left(2 x^1+x^2\right) \left(x^3+x^4\right)}{4 (x^1)^2 x^3 (x^2)^3+4 (x^1)^2 x^4 (x^2)^3-4 \left(x^1+x^2\right) (x^3)^2 (x^4)^2 x^2}&\mbox{for} \ i = 2\\

&\frac{\left(x^1+x^2\right) \left(x^3+2 x^4\right)}{4 x^3 \left(-(x^1)^2 \left(x^3+x^4\right) (x^2)^2+(x^3)^2 (x^4)^2 x^2+x^1 (x^3)^2 (x^4)^2\right)}&\mbox{for} \ i = 3\\

&\frac{\left(x^1+x^2\right) \left(2 x^3+x^4\right)}{4 x^4 \left(-(x^1)^2 \left(x^3+x^4\right) (x^2)^2+(x^3)^2 (x^4)^2 x^2+x^1 (x^3)^2 (x^4)^2\right)}&\mbox{for} \ i = 4.\\
\end{array}\right.
\ee
Now one can calculate exterior derivatives and wedge products between $\Pi$, $\Phi$, $\Psi$ and $\Theta$, and then easily check from \eqref{eq3.5} that the manifold satisfies $R\cdot R =0$, i.e., semisymmetric, although the 1-forms are not closed, since $\Theta$ being arbitrary.\\
Again as the manifold is $S$-$GK_n$ and also $SGK_n$, by Theorem \ref{th3.8}, it is $C$-$SGK_n$, $W$-$SGK_n$ and also $K$-$SGK_n$ but of different 1-forms, which follows from Theorem \ref{th3.8} and Corollary \ref{cor3.1} as $M$ is not Ricci recurrent.
\section{\bf{Conclusion}}
In this paper we study the geometric properties of a super generalized recurrent manifold. It is shown that its associated 1-forms are not uniquely determined and they are linearly dependent with $d\kappa$ and also their dependency relations are evaluated. We found out the form of $R\cdot R$ of a $SGK_n$ and showed that it is semisymmetric if all of its associated 1-forms are closed and pairwise codirectional. It is also shown that if the manifold is Roter type then super generalized recurrent and Ricci generalized recurrent manifolds are equivalent. Again we prove that Ricci recurrency is a necessary and sufficient condition for a $SGK_n$ to be $C$-$SGK_n$ or $P$-$SGK_n$ or $K$-$SGK_n$ with same associated 1-forms. Finally a proper example of $SGK_n$ is given which verifies the main results of the paper.
\newline\noindent
\textbf{Acknowledgment:} 
The author HK gratefully acknowledges to CSIR, New Delhi (File No. 09/025 (0194)/2010-EMR-I) for the financial assistance. All the algebraic computations of Section \ref{exam} are performed by a program in Wolfram Mathematica.


\end{document}